\documentclass[11pt]{amsart}
\usepackage{amsmath, amssymb, array}
\usepackage[top=30truemm,bottom=30truemm,left=30truemm,right=30truemm]{geometry} 
\usepackage{mathtools}
\usepackage[all]{xy}
\usepackage{ascmac}
\usepackage[dvips]{graphicx}
\usepackage{color}
\usepackage{amscd}
\usepackage{fancyhdr}
\usepackage[dvipdfmx]{hyperref}
\usepackage{amsrefs}
\usepackage{cleveref}
\usepackage{chngcntr}
\usepackage{apptools}
\usepackage{comment}

\newcommand{\fn}{\footnote}

\renewcommand{\Re}{\mathrm{Re}} %real part

\DeclareMathOperator{\im}{Im} %image/imaginary part
\DeclareMathOperator{\rank}{rank} %rank
\title{Quadrature formulas for Bessel polynomials}

\author{Hideki Matsumura$^{*}$}
\email{hidekimatsumura@keio.jp}
\address{${}^*$Department of Mathematics, Faculty of Science and Technology, Keio University, 3-14-1 Hiyoshi , Kouhoku-ku, Yokohama, Kanagawa, Japan}

\thanks{This research is supported by KAKENHI 18H05233.}
 
\subjclass[2010]{primary 33C45; secondary 65D32; tertiary 14G05}
\keywords{quadrature formula, Bessel polynomial, Riesz--Shohat theorem, Newton polygon, Christoffel--Darboux kernel, 
rational points on elliptic curves}
\date{\today}

\theoremstyle{plain}
 \newtheorem{theorem}{Theorem}[section] 
  \crefname{theorem}{Theorem}{Theorems}
 \newtheorem{proposition}[theorem]{Proposition}
 \crefname{proposition}{Proposition}{Propositions}
 \newtheorem{lemma}[theorem]{Lemma}
 \crefname{lemma}{Lemma}{Lemmas}
 
  \crefname{corollary}{Corollary}{Corollaries}
 
   \crefname{conjecture}{Conjecture}{Conjectures}
 
 \crefname{question}{Question}{Questions}
 \newtheorem{problem}[theorem]{Problem}
   \crefname{problem}{Problem}{Problems}
 
    \crefname{notation}{Notation}{Notations}
 
\theoremstyle{definition} 
 \newtheorem{definition}[theorem]{Definition}
  \crefname{definition}{Definition}{Definitions}
 \newtheorem{example}[theorem]{Example}
   \crefname{example}{Example}{Examples}
 \newtheorem{remark}[theorem]{Remark}
   \crefname{remark}{Remark}{Remarks}
    
   \crefname{claim}{Claim}{Claims}

\begin{document}

%%%Title&Author

\maketitle

\tableofcontents
%%% Abstract 

\begin{abstract}
A quadrature formula is a formula computing a definite integration by evaluation at finite points. 
The existence of certain quadrature formulas for orthogonal polynomials is related to interesting problems such as Waring's problem in number theory and spherical designs in algebraic combinatorics.
Sawa and Uchida proved the existence and the non-existence of certain rational quadrature formulas for the weight functions of certain classical orthogonal polynomials. 
Classical orthogonal polynomials belong to the Askey-scheme, which is a hierarchy of hypergeometric orthogonal polynomials. 
Thus, it is natural to extend the work of Sawa and Uchida to other polynomials in the Askey-scheme.
In this article, we extend the work of Sawa and Uchida to the weight function of the Bessel polynomials. 
In the proofs, we use the Riesz--Shohat theorem and Newton polygons.
It is also of number theoretic interest that proofs of some results are reduced to determining the sets of rational points on elliptic curves.
\end{abstract}

\section{Introduction}
Let $\gamma \subset \mathbb{C}$ be a smooth path and $x_1,\ldots, x_m$, $y_1,\ldots, y_m \in \mathbb{C}$.
Suppose that $w(z)$ is a function such that $\int_{\gamma} f(z) w(z)dz$ exists for any $f(z) \in \mathbb{C}[z]$.
An integration formula of the form
\[\sum_{i=1}^{m} x_i f(y_i)=\int_{\gamma} f(z) w(z)dz\]
is called a quadrature formula of degree $n$ for a weight function $w(z)$ if it holds for all $f(z) \in \mathbb{C}[z]$ such that $\deg(z) \leq n$.
 The points $y_i$'s are called the nodes of the quadrature formula.
 For example, for a real interval $\gamma=[a,b]$, 
  $w(t)=1$ and $n=1$, the trapezoidal rule states that
 \[\frac{b-a}{2}(f(a)+f(b))=\int_a^b f(t) dt \]
is a quadrature formula of degree $1$ for $w(t)=1$ with $2$ nodes.
 For $\gamma=[a,b]$, $w(t)=1$ and $n=3$, the Simpson rule states that
 \[\frac{b-a}{6} \left(f(a)+4 f\left(\frac{a+b}{2} \right)+f(b) \right)=\int_a^b f(t) dt\]
 is a quadrature formula of degree $3$ for $w(t)=1$ with $3$ nodes.
 
A quadrature formula is called rational if all nodes are rational numbers. The existence of certain rational quadrature formulas for orthogonal polynomials is related to interesting problems such as Waring's problem in number theory (cf.\ \cite{Hausdorff,Nesterenko}) and spherical designs in algebraic combinatorics (cf.\ \cite{BB2005,DGS1977}). 
See \cite{Hausdorff,SU20} for applications.

Quadrature formulas are determined by the weight function $w(t)$, the degree $n$ and the number $m$ of nodes.
The following problem is fundamental:

\begin{problem}
For a fixed weight function $w(t)$ and a natural number $m$, determine the range of $n$ such that
there exists a rational quadrature formula of degree $n$ for $w(t)$ with $m$ nodes.
\end{problem}

By the Stroud-type bound (cf.\ \cite[Proposition 4.3]{SU19}, \cite[Proposition 2.7]{SU20}), if there exists a quadrature formula of degree $n$ for $w(t)$ with $m$ nodes, then $n \leq 2m-1$. 
The case $n=2m-1$ is called the tight case.
Sawa and Uchida \cite{SU19,SU20} considered the ``almost tight case" and proved that there exists no quadrature formula of degree $2r$ for certain classical orthogonal polynomials (Hermite polynomials $H_n(t)$, Legendre polynomials $P_n(t)$ and Laguerre polynomials $L_n^{(0)}(t)$) with $r+1$ nodes on $\mathbb{Q}$ (\cite[Theorem 1.4]{SU20}).
 They also proved that there exist quadrature formulas of degree $r+1$ for the above classical orthogonal polynomials with $r+1$ nodes on $\mathbb{Q}$ (\cite[Theorem 5.5]{SU20}). 
 They took algebro-geometric approaches involving the Riesz--Shohat theorem, Christoffel--Darboux kernels and Newton polygons.

Classical orthogonal polynomials belong to the Askey-scheme \cite{KLS}, which is a hierarchy of hypergeometric orthogonal polynomials.
Thus, it is natural to extend the work of \cite{SU19,SU20} to more general polynomials in the Askey-scheme.
In this article, we extend the work of Sawa and Uchida \cite{SU19,SU20} for the weight function of Bessel polynomials \cite[p.\ 101]{KF}
\begin{align*}
y_n(x):=\sum_{k=0}^n \frac{(n+k)!}{(n-k)!k!} \left(\frac{x}{2} \right)^k.
\end{align*}
The Bessel polynomials belong to the Askey-scheme and
are orthogonal with respect to the weight function 
\[w(z)=-\frac{1}{4 \pi \sqrt{-1}} e^{-\frac{2}{z}},\]
 where the path of integration is the unit circle $S^1 \subset \mathbb{C}$, i.e.,
 \[\int_{S^1} y_m(z)y_n(z)w(z)dz=0\]
 for $m \neq n$ (cf.\ \cite[p.\ 104]{KF}).

Orthogonal polynomials are defined from certain inner products expressed as integrations.
The inner products for the classical orthogonal polynomials in \cite{SU19,SU20} are expressed as real integrations. 
On the other hand, in our case, note that the inner product is expressed as complex integrations.
As a generalization of  ``quadrature formulas with nodes on $\mathbb{Q}$" (\cite{SU19,SU20}), 
we consider ``quadrature formulas with the real parts and the imaginary parts of nodes on $\mathbb{Q}$."
Since the path of integration is $S^1$, we consider the following problems:

\begin{problem} \label{problem} 
\begin{enumerate}
\item Does there exist a quadrature formula of degree $2r$ for Bessel polynomials with $r+1$ nodes on $\mathbb{Q}(\sqrt{-1}) \cap S^1$?
\item Does there exist a quadrature formula of degree $r+1$ for Bessel polynomials with $r+1$ nodes on $\mathbb{Q}(\sqrt{-1}) \cap S^1$?
\end{enumerate}
\end{problem}

There are some previous works concerning spherical designs and quadrature formulas on large enough extensions of $\mathbb{Q}$.
For example, for fixed $t$ and $d \in \mathbb{Z}_{>0}$, Cui, Xia and Xiang \cite[Theorem 1.3]{CXX} proved the existence of spherical $t$-designs on $S^{d-1}$
with large enough nodes on the algebraic extension $\mathbb{Q}(\sqrt{q} \mid q: \mbox{prime})$ over $\mathbb{Q}$.
For fixed $p \in \mathbb{Z}_{\geq 2}$ and $n \in \mathbb{Z}_{>0}$, Kuperberg \cite[Theorem 1, p.\ 855]{Kuperberg} proved the existence of certain 
Chebyshev-type quadrature formulas (quadrature formulas with the same weights) of degree $2n+1$ with $p^n$ nodes on certain number fields 
(see also \cite[Lemma 2]{Kuperberg}).
As far as the author knows, there are few works considering quadrature formulas whose nodes lie in a quadratic extension of $\mathbb{Q}$.
Note that \cite{SU19,SU20} do not consider quadrature formulas with nodes on $S^1$. 
In this article, we give negative answers to \cref{problem} (1) for all $r \in \mathbb{Z}_{\geq 1}$ and \cref{problem} (2) for $r=2$.

\begin{theorem} \label{MT1} 
\begin{enumerate} 
\item There exists no quadrature formula of degree $2$ for Bessel polynomials with $2$ nodes on $\mathbb{Q}(\sqrt{-1}) \cap S^1$.
\item There exists no quadrature formula of degree $4$ for Bessel polynomials with $3$ nodes on $\mathbb{Q}(\sqrt{-1}) \cap S^1$.
\item 
There exists no quadrature formula of degree $2r$ for Bessel polynomials with $r+1$ nodes on $\mathbb{Q}(\sqrt{-1})$ for all $r \in \mathbb{Z}_{\geq 3}$.
\end{enumerate}
\end{theorem}

The proof of (1) is straightforward.
(2) follows from \cref{MT2}.
In the proof of (3), we use the Riesz--Shohat theorem (\cref{RS}) and a fact on Newton polygons (\cref{NP}).
An advantage of Bessel polynomials is that the computation (or estimate) of $p$-adic valuation is relatively easy. 
On the other hand, to prove (2), we cannot use \cref{NP} since  
we can construct infinitely many quadrature formulas of degree $4$ for Bessel polynomials with $3$ nodes on $\mathbb{Q}(\sqrt{-1})$
by using the Christoffel--Darboux kernel (cf.\ \cref{4-3}).
By contrast, there exists no quadrature formula of degree $4$ for Hermite polynomials with $3$ nodes on $\mathbb{Q}$ by \cite[Theorem 1.4]{SU20}. 

\begin{theorem} \label{MT2}
There exists no quadrature formula of degree $3$ for Bessel polynomials with $3$ nodes on $\mathbb{Q}(\sqrt{-1}) \cap S^1$.
\end{theorem}
Note that it is an immediate consequence of the Riesz--Shohat theorem that
there exists a quadrature formula of degree $r$ for Bessel polynomials with $r+1$ nodes on $\mathbb{Q}(\sqrt{-1}) \cap S^1$ (\cref{RSrem}). 
Thus, increasing the degree is a problem.
In the proof of \cref{MT2}, we use the Riesz--Shohat theorem (\cref{RS}) 
to reduce the problem to determining the set of rational points on an elliptic curve.
It seems hard to extend the proof of \cref{MT2} to $r \geq 3$.
On the other hand, by the Riesz--Shohat theorem (\cref{RS}),
we can prove that there exist quadrature formulas of degree $r+1$ with $r+1$ nodes on $\mathbb{Q}$ for all $r \in \mathbb{Z}_{\geq 0}$.
\begin{theorem} \label{MT3}
There exist quadrature formulas of degree $r+1$ for Bessel polynomials with $r+1$ nodes on $\mathbb{Q}$ for all $r \in \mathbb{Z}_{\geq 0}$.
\end{theorem}

This article is organized as follows:

In \S2, we formulate a quadrature formula on $\mathbb{C}$ and introduce three key tools: 
the Riesz--Shohat theorem (\cref{RS}), Newton polygons and Christoffel--Darboux kernels.
In \S3, we prove \cref{MT1} by using the Riesz--Shohat theorem and Newton polygons.
We also construct an example of quadrature formulas of degree $4$ for Bessel polynomials with $3$ nodes on $\mathbb{Q}(\sqrt{-1})$ 
by using the Christoffel--Darboux kernel (\cref{4-3}).
In \S4, we prove \cref{MT2} by using the Riesz--Shohat theorem. 
It is also of number theoretic interest that this proof is reduced to determining the set of rational points on an elliptic curve.
Then, in \S5, we prove
\cref{MT3} by using the Riesz--Shohat theorem.

\section{Quadrature formulas for Bessel polynomials}

We consider the ``quadrature formulas" for Bessel polynomials. Since
\[\frac{1}{2^k} \frac{(n+k)!}{(n-k)!k!}=\frac{1}{2^k} \frac{(n+k)!}{(n-k)!(2k)!} \frac{(2k)!}{k!}=\binom{n+k}{2k}(2k-1)!!,\]
note that $y_n(x) \in \mathbb{Z}[x]$.

Unlike the classical orthogonal polynomials in \cite{SU20},
the path of integration is $S^1 \subset \mathbb{C}$.
Let $\gamma \subset \mathbb{C}$ be a smooth path. %fn
\fn{Later, we consider $\gamma=S^1$. 
For the formulation over $\mathbb{R}$, see \cite{SU20}.
} %end
Let $w(z)$ be a complex function such that
$\int_\gamma f(z) w(z) dz$ exists for any $f(z) \in \mathbb{C}[z]$. Suppose that 
$\int_\gamma \phi_l(z)^2 w(z) dz \neq 0$ for all $l$, where $\{\phi_l\}$ is the system of monic orthogonal polynomials
with respect to $w(z)$.
We consider ``quadrature formulas" in the following sense:

\begin{definition} \label{CQF}
An integration formula of the form
\begin{align} \label{QF}
\sum_{i=1}^{m} x_i f(z_i)=\int_{\gamma} f(z) w(z)dz
\end{align}
with $x_1,\ldots, x_m, z_1, \ldots, z_m \in \mathbb{C}$
is called a quadrature formula of degree $n$ for a weight function $w(z)$ with nodes $z_i$ and weights $x_i$ if 
it holds for all $f(z) \in \mathbb{C}[z]$ such that $\deg(f) \leq n$.
A quadrature formula on a subset $S \subset \mathbb{C}$ is a quadrature formula such that $z_i \in S$ for all $i$. 
\end{definition}

For Bessel polynomials, we take 
\[w(z)=-\frac{1}{4 \pi \sqrt{-1}} e^{-\frac{2}{z}} \]
so that 
\[  \int_{S^1} w(z) dz=1. 
\] 

We introduce key tools: 
the Riesz--Shohat theorem \cite{Shohat} and a lemma on Newton polygons.
They are key tools in the proof of \cite[Theorem 1.4, Theorem 5.5]{SU20}, and also work over $\mathbb{C}$. 

In the proof of \cref{MT1,MT2,MT3}, we use the Riesz--Shohat theorem.
The proof is exactly similar to the real case (\cite{Shohat}).
\begin{theorem} [{Riesz--Shohat theorem, cf.\ \cite[Proposition 2.4]{SU20}, \cite[Theorem I]{Shohat}}] \label{RS} %[SU20] p. 1246, 1247, [Shohat] p. 465
Suppose that $1 \leq k \leq r+2$, $z_1, \ldots, z_{r+1} \in \mathbb{C}$ are distinct and let
\[\theta_{r+1}(z) :=\prod_{i=1}^{r+1}(z-z_i).\]
Then, the following are equivalent:
\begin{enumerate}
\item There exist $x_1, \ldots, x_{r+1} \in \mathbb{C}$ such that
\[\sum_{i=1}^{r+1} x_i f(z_i)=\int_{\gamma} f(z)w(z)dz\]
is a quadrature formula of degree $2(r+1)-k$.
\item $\theta_{r+1}$ is a quasi-orthogonal polynomial of degree $r+1$ and of order $k-1$, i.e.,
there exist $b_1, \ldots, b_{k-1} \in \mathbb{C}$ such that
\[\theta_{r+1}(z)=\phi_{r+1}(z)+b_1 \phi_r(z)+ \cdots +b_{k-1} \phi_{r+2-k}(z).\]
Here, $\{ \phi_l \}$ is the system of monic orthogonal polynomials with respect to $w(z)$.
\end{enumerate}
Furthermore, if the above equivalent conditions hold, then we have
\[x_i =\int_{\gamma} \frac{\theta_{r+1}(z)}{(z-z_i)\theta'_{r+1}(z_i)}w(z)dz.\]
\end{theorem}
In what follows, $\{\phi_l \}$ denotes the system of monic orthogonal polynomials with respect to a weight function $w(t)$ and
$\{\Phi_l \}$ denotes a system of (not necessarily monic) orthogonal polynomials with respect to $w(t)$.

\begin{remark} \label{RSrem}
Let $w(z)$ be a weight function and $A$ be a subset of $\mathbb{C}$ with at least $r+1$ elements.
Then, there exists a quadrature formula of degree $r$ for $w(z)$ with $r+1$ nodes on $A$. 
Indeed, this is an immediate consequence of \cref{RS}. 
Note that this is the case $2r+2-k=r$. i.e., $k=r+2$ in \cref{RS}. 
Let $\{\phi_l\}$ be the system of monic orthogonal polynomials with respect to $w(z)$. 
Then, since $\{\phi_0, \ldots, \phi_{r+1}\}$ is a basis of the space of polynomials of degree $\leq r+1$
and $\theta_{r+1}(z)$ is monic,
there exist $b_1, \ldots, b_{r+1} \in \mathbb{C}$ such that
\[\theta_{r+1}(z)=\phi_{r+1}(z)+b_1 \phi_r(z)+ \cdots +b_{r+1} \phi_0(z).\]
Therefore, by \cref{RS}, for all $z_1, \ldots, z_{r+1} \in A$,
there exist $x_1, \ldots, x_{r+1} \in \mathbb{C}$ such that
\[\sum_{i=1}^{r+1} x_i f(z_i)=\int_{\gamma} f(z)w(z)dz\]
is a quadrature formula of degree $r$.
In particular, there exists a quadrature formula of degree $r$ for Bessel polynomials with $r+1$ nodes on $\mathbb{Q}(\sqrt{-1}) \cap S^1$. %fn
\fn{We can take arbitrary $r+1$ points on $\mathbb{Q}(\sqrt{-1}) \cap S^1$ as nodes.} %end
\end{remark}

In the proof of \cref{MT1} (3), we use a fact on Newton polygons.

\begin{lemma} [{Cf.\ \cite[Lemma 4.2]{SU20}}] \label{NP} %\label{KL} p. 1253 
Let $K$ be a discrete valuation field, $\mathcal{O}_K$ be its ring of integers and $f(x) \in \mathcal{O}_K[x]$. If all zeros of $f(x)$ are $K$-rational, 
then all edges of the Newton polygon of $f(x)$ 
have integral slopes. 
\end{lemma}

We also use the following version of Bertrand's postulate:

\begin{theorem} [{\cite[p.\ 505]{Breusch}}] \label{BPcong} 
For $n \geq 4$,
there exists a prime number $p$ such that $n<p \leq 2n$ and $p \equiv 3 \pmod{4}$. 
\end{theorem}

As in the proof of \cite[Lemma 4.5]{SU20}, we obtain the following lemma:

\begin{lemma} \label{BPC}
For all $l \in \mathbb{Z}_{\geq 7}$, there exists a prime number $p$ such that $(l+1)/2<p \leq l$ and $p \equiv 3 \pmod{4}$.
\end{lemma}

In \cref{4-3}, we use the Christoffel--Darboux kernel.

\begin{definition} [{\cite[p.\ 1243]{SU20}}] 
Let $\{\Phi_l \}$ be a system of orthogonal polynomials with respect to a weight function $w(z)$ over $\gamma \subset \mathbb{C}$.
The Christoffel--Darboux kernel for polynomials of degree at most $l$ is defined by
\[K_l(x,y):=\sum_{k=0}^l h_lh_k^{-1} \Phi_k(x) \Phi_k(y) .\]
Here, $h_k:=\int_{\gamma} \Phi_k(z)^2 w(z) dz$.
\end{definition}

The Christoffel--Darboux kernel is computed by the following formula:

\begin{proposition} [{Chirstoffel--Darboux formula, cf.\ \cite[Proposition 2.1]{SU20}, \cite[Theorem 3.2.2]{Szego}}] \label{CDF} %[Szego] p. 43
For $l \in \mathbb{Z}_{\geq 0}$, let $k_l$ be the leading coefficient of $\Phi_l(x)$. Then,
\[K_l(x,y)=\frac{k_l}{k_{l+1}} \frac{\Phi_{l+1}(x) \Phi_l(y)-\Phi_l(x) \Phi_{l+1}(y)}{x-y}. \]
\end{proposition}

\begin{definition} \label{fr}
\[f_l(x,y):=\frac{\Phi_{l+1}(x) \Phi_l(y)-\Phi_l(x) \Phi_{l+1}(y)}{x-y}.\]
\end{definition}
By \cref{RS} and \cref{CDF}, we obtain the following lemma:

\begin{lemma} [{Cf.\ \cite[Lemma 2.5]{SU20})}] \label{CDK}
Suppose that $x_1, \ldots, x_{r+1}$, $z_1, \ldots , z_{r+1} \in \mathbb{C}$ satisfy 
\cref{QF} for $(m,n)=(r+1,2r)$.
Assume that $z_1, \ldots , z_{r+1}$  are distinct. 
Then $f_r(z_i,z_j) = 0$ for every distinct $i$, $j$.
\end{lemma}

\section{Proof of \cref{MT1}}
In this section, we prove \cref{MT1}. 
We also construct an example of quadrature formulas of degree $4$ for Bessel polynomials with $3$ nodes on $\mathbb{Q}(\sqrt{-1})$.

\subsection{Proof of \cref{MT1} (1), (2)}
In this subsection, we prove \cref{MT1} (1) and (2).

\begin{proof} [{Proof of \cref{MT1} (1), (2)}]
\begin{enumerate}
\item Suppose that there exist 
$x_1$, $x_2 \in \mathbb{C}$, $z_1$, $z_2 \in \mathbb{Q}(\sqrt{-1})$ such that $|z_1|=|z_2|=1$ and
\begin{align*} 
x_1+x_2 &=1, \\
x_1z_1+x_2z_2 &=-1, \\
x_1z_1^2+x_2z_2^2 &=\frac{2}{3}.
\end{align*}
By the first and the second equality, we obtain 
\[(z_1-z_2)x_1= -z_2-1.\]
If $z_1=z_2$, then 
\[z_1=z_2=-1.\]
Thus, by the third equality,
\[x_1+x_2=\frac{2}{3},\]
 which contradicts the first equality. 
Since $z_1 \neq z_2$,
 \[x_1=\frac{-z_2-1}{z_1-z_2}.\]
 We substitute it and the first equality into the third equality to obtain
 \[ -z_1-z_1z_2-z_2=\frac{2}{3}. \] 
If $z_1=-1$, then
 \[1+z_2-z_2=\frac{2}{3}, \]
 which is a contradiction. Thus, $z_1 \neq -1$.
 Therefore, 
\[z_2=\frac{-z_1-\frac{2}{3}}{z_1+1}. \]
Since $|z_1|=|z_2|=1$,
\[ \Re(z_1)=-\frac{5}{6}.\]  
Thus, we have
\[\im(z_1)= \pm \frac{\sqrt{11}}{6} \not\in \mathbb{Q}, \]
which is a contradiction.

\item This is an immediate consequence of \cref{MT2}. 
\end{enumerate}
\end{proof}

\begin{remark} \label{UQ}
There exists a unique quadrature formula of degree $2$ for Bessel polynomials with $2$ nodes on $S^1 \subset \mathbb{C}$.
 Indeed, if there exists such a quadrature formula on $S^1$, then we have 
 \begin{align*}
 z_1 &=\frac{-5+\sqrt{-11}}{6},\\
 z_2 &=\frac{-5-\sqrt{-11}}{6}
 \end{align*}
 by the proof of \cref{MT1} (1).
Then, by \cref{RS} and the residue theorem, 
\begin{align*}
x_1 &=\int_{S^1} \frac{z-z_2}{z_1-z_2} w(z) dz=\frac{11+\sqrt{-11}}{22}, \\
x_2 &=\int_{S^1} \frac{z-z_1}{z_2-z_1} w(z) dz=\frac{11-\sqrt{-11}}{22}.
\end{align*}
Since
\[-\frac{1}{4 \pi \sqrt{-1}}\int_{S^1} z^j e^{-\frac{2}{z}} dz=\frac{(-2)^j}{(j+1)!}\]
for $j \in \mathbb{Z}_{\geq 0}$ by the residue theorem, we can check that 
\[x_1f(z_1)+x_2f(z_2)=-\frac{1}{4 \pi \sqrt{-1}}\int_{S^1} f(z) e^{-\frac{2}{z}} dz\] 
is a quadrature formula of degree $2$ for Bessel polynomials.
\end{remark}

\subsection{Proof of \cref{MT1} (3)}
In this subsection, we prove \cref{MT1} (3). 
We also construct an example of quadrature formulas of degree $4$ for Bessel polynomials with $3$ nodes on $\mathbb{Q}(\sqrt{-1})$.

\begin{proof} [{Proof of \cref{MT1} (3)}]

Assume that there exist $x_1, \ldots, x_{r+1} \in \mathbb{C}, z_1, \ldots, z_{r+1} \in \mathbb{Q}(\sqrt{-1})$ such that
\[\sum_{i=1}^{r+1} x_iz_i^j= -\frac{1}{4 \pi \sqrt{-1}}  \int_{S^1} z^j e^{-\frac{2}{z}}  dz \quad (j=0,1, \ldots, 2r). 
\]
By \cref{RS}, there exists $c \in \mathbb{Q}(\sqrt{-1})$ such that the zeros of 
the quasi-Bessel polynomial $y_{r+1;c}(x)$ are $z_1, \ldots, z_{r+1}$.
Write $c=s/t$ with $s, t \in \mathbb{Z}[\sqrt{-1}]$ and $\gcd (s,t)=1$. Let
\[f(x):=ty_{r+1;c}(x)= ty_{r+1}(x)+sy_{r}(x) \in \mathbb{Z}[x].\]
Write
\[f(x)=\sum_{k=0}^{r+1} a_k x^{r+1-k}. \]
We have
\begin{align} \label{(25)'} 
a_k=\begin{cases}
t(2r+1)!! & (k=0),\\
(t(2r+1) +s)(2r-1)!! & (k=1), \\
\left(t \binom{2(r+1)-k}{2(r+1-k)}+s \binom{2r+1-k}{2(r+1-k)} \right) (2(r+1-k)-1)!! & (2 \leq k \leq r+1). 
\end{cases}
\end{align}
First, assume that $r \geq 7$, $r=3$ or $4$. 
For $r \geq 7$, by \cref{BPC}, we can take a prime number $p$ such that $(r+1)/2<p \leq r$ and $p \equiv 3 \pmod{4}$ (i.e., $p$ is inert in $\mathbb{Q}(\sqrt{-1})$). 
For $r=3$, $4$, let $p=3$.
We consider the Newton polygon of $f(x)$ with respect to $p$.
Let $P_k:=(k,v_p(a_k)) \in \mathbb{R}^2$.

Let
\[j:=\begin{cases}
r & (t+s \equiv 0 \pmod{p}),\\
r+1 & (t+s \not\equiv 0 \pmod{p}).
\end{cases}\]
Then, $v_p(a_j)=0$ by \cref{(25)'}.
Indeed, if $t+s \equiv 0 \pmod{p}$, then $v_p(s)=v_p(t)=0$ by $\gcd (s,t)=1$.
Since $p<r+1<2p$ and $p$ is odd, 
\[v_p(a_r)=v_p \left(\frac{(r+1)((t+s)r+2t)}{2} \right)=0.\]
If $t+s \not\equiv 0 \pmod{p}$, then $v_p(a_{r+1})=v_p(t+s)=0$.

Therefore, there exists $j_0:=\min \{j \mid v_p(a_j)=0 \}$. 

If $v_p(t) =0$, then 
\[v_p(a_0)=v_p((2r+1)!!)= \begin {cases}
1 & (2r+1<3p),\\
2 & (2r+1 \geq 3p, \; p>3),\\ 
3 & (2r+1 \geq 3p, \; p=3) 
\end{cases}\] 
by \cref{(25)'} and $p<2r+1<4p$.
If $v_p(t) \geq 1$, then $v_p(s)=0$ by $\gcd (s,t)=1$. Thus, %fn
\fn{If $p=3$, then $r<5$ by $(r+1)/2<p$. Thus, $2r-1<9=3p$.} %end
\[v_p(a_1)=v_p((2r-1)!!)= \begin {cases}
1 & (2r-1<3p), %\fnm
\\ 
2 & (2r-1 \geq 3p) 
\end{cases}\] %fn
%\fnt{If $p=3$, then $r<5$ by $(r+1)/2<p$. Thus, $2r-1<9=3p$.} %end
 by \cref{(25)'} and $p<2r-1<4p$.

\begin{enumerate}
\item
Suppose that $v_p(t)=0$.

\begin{enumerate}
\item
 If $2r+1<3p$, %fn
 \fn{This includes the case $r=3$.} %end
 then $v_p(a_0)=1$. Let $j_0:=\min \{j \mid v_p(a_j)=0 \}$. 
In this case, the Newton polygon of $f(x)$ has the edge $P_0P_{j_0}$, 
whose slope is $-1/j_0 \not\in \mathbb{Z}$
since $j_0 > r-(p-1)/2 \geq p-(p-1)/2 =(p+1)/2 \geq 2$ by $p \geq 3$.
For the first inequality, note that if $p$ is an odd prime and $p \leq 2(r+1-k)-1$, i.e., $k \leq r-(p-1)/2$, then $v_p(a_k) \geq 1$.

\item
If $r \geq 7$ and $2r+1 \geq 3p$, then  $v_p(a_0)= 2$ and $v_p(a_1) \geq 1$. 
\begin{enumerate}

\item
If $v_p(a_1)=1$, then the Newton polygon of $f(x)$
has the edge $P_1P_{j_0}$, whose slope is $-1/(j_0-1) \not\in \mathbb{Z}$ since $j_0 > 2$.

\item If $v_p(a_1) \geq 2$ 
and $v(a_i) \neq 1$ for all $i$, then
the Newton polygon of $f(x)$
has the edge $P_0P_{j_0}$, whose slope is $-2/j_0 \not\in \mathbb{Z}$ since $j_0 > 2$. 

\item If $v_p(a_1) \geq 2$ 
and $v(a_i)=1$ for some $i \geq 2$, then let 
$i_0:=\min \{i \mid v_p(a_i)=1 \}$. 
The Newton polygon of $f(x)$
has the edge $P_0P_{i_0}$, whose slope is $-1/i_0 \not\in \mathbb{Z}$ since $i_0 >1$.
\end{enumerate}

\item If $r=4$, then
\begin{align*} 
a_0 &=945t, \\
a_1 &=3(315t+35s), \\
a_2 &= 3(140t+35s), \\
a_3 &=3(35t+15s), \\
a_4 &=15t+10s, \\
a_5 &=t+s.
\end{align*}

\begin{enumerate}
\item If $v_3(s)>0$, then $v_3(a_0)=3$, 
$v_3(a_1) \geq 2$, $v_3(a_2)=1$, $v_3(a_3)$, $v_3(a_4) \geq 1$ and $v_3(a_5)=0$.
Thus, the Newton polygon of $f(x)$ has the edge $P_2P_5$, 
whose slope is $-1/3 \not\in \mathbb{Z}$.

\item If $v_3(s)=0$, then 
$v_3(a_0) = 3$, $v_3(a_1)=1$, $v_3(a_2)$, $v_3(a_3) \geq 1$ and $v_3(a_4)=0$.
Thus, the Newton polygon of $f(x)$ has the edge $P_1P_4$, 
whose slope is $-1/3 \not\in \mathbb{Z}$.
\end{enumerate}

\end{enumerate}

\item
Suppose that $v_p(t)>0$. Since $\gcd (s,t)=1$, $v_p(s)=v_p(t+s)=0$.
We have $v_p(a_0) \geq 2$.

\begin{enumerate}
\item
If $2r-1<3p$, then $v_p(a_1)=1$ and the Newton polygon of $f(x)$
has the edge $P_1P_{j_0}$, whose slope is $-1/(j_0-1) \not\in \mathbb{Z}$
since $j_0 > 2$.

\item
 If $2r-1 \geq 3p$, then $v_p(a_0) \geq 3$ and $v_p(a_1)=2$.
We have $v_p(a_2) \geq 1$ by $j_0>(p+1)/2 \geq 3$.
 \begin{enumerate}
\item
If $v_p(a_2)=1$, then the Newton polygon of $f(x)$
has the edge $P_2P_{j_0}$, whose slope is $-1/(j_0-2) \not\in \mathbb{Z}$ 
since $j_0 > r-(p-1)/2 \geq (3p+1)/2-(p-1)/2 \geq p+1 \geq 4$ by $2r-1 \geq 3p$.
The first inequality follows from the argument in Case (1-a).

\item If $v_p(a_2) \geq 2$ and $v_p(a_i) \neq 1$ for all $i$, then
the Newton polygon of $f(x)$
has the edge $P_1P_{j_0}$, whose slope is $-2/(j_0-1) \not\in \mathbb{Z}$ since $j_0 > 4$.

\item If $v_p(a_2) \geq 2$ and $v_p(a_i) = 1$ for some $i \geq 3$, then 
the Newton polygon of $f(x)$ has the edge $P_1P_{i_0}$, 
whose slope is $-1/(i_0-1) \not\in \mathbb{Z}$ since $i_0 >2$. 
\end{enumerate}
\end{enumerate}

\end{enumerate}

For $r=5$, $6$, we consider $(2+\sqrt{-1})$-adic valuation and we can prove the theorem similarly. %fn
\fn{Note that for $n \in \mathbb{Z}$, $n$ is divisible by $2+\sqrt{-1}$ in $\mathbb{Z}[\sqrt{-1}]$ if and only if $n$ is divisible by $5$ in $\mathbb{Z}$.} %end
\end{proof}

Next, we construct an example of quadrature formulas of degree $4$ for Bessel polynomials with $3$ nodes on $\mathbb{Q}( \sqrt{-1})$.

Since $y_2(z)=3z^2+3z+1$ and $y_3(z)=15z^3+15z^2+6z+1$,
\begin{align*}
f_2(x,y) = 3((15y^2+15y+5)x^2+(15y^2+14y+4)x+(5y^2+4y+1)). 
\end{align*}
by \cref{fr}.
Multiplying both sides of $f_2(x,y)=0$ by $4(15y^2+15y+5)$, we have
\[(2(15y^2+15y+5)x+(15y^2+14y+4))^2+75y^4+120y^3+84y^2+28y+4=0. \]
Let $w:=2(15y^2+15y+5)x+(15y^2+14y+4)$. Then we have a curve
\[C: w^2=-75y^4-120y^3-84y^2-28y-4. \]
The problem is reduced to determining the set $C(\mathbb{Q}(\sqrt{-1}))$
of $\mathbb{Q}(\sqrt{-1})$-rational points on $C$.
Let $J$ be the Jacobian variety of $C$. Then,
by MAGMA \cite{Bosma-Cannon-Playoust}, we can check that $\rank J(\mathbb{Q}(\sqrt{-1}))=1$.

\begin{example} \label{4-3}
From the rational points on the above elliptic curve $C$, we can construct an example of quadrature formulas of degree $4$ for Bessel polynomials with $3$ nodes on $\mathbb{Q}(\sqrt{-1})$. 
For example, by MAGMA \cite{Bosma-Cannon-Playoust}, 
\[P:=(-264/743, \pm 377866 \sqrt{-1}/552049) \in C(\mathbb{Q}(\sqrt{-1})).\]
Then, by Mathematica \cite{Wolfram},
for any weight $x_1$, other weights $x_2$, $x_3$ and nodes $z_1$, $z_2$, $z_3$ are expressed as rational functions of $x_1$ and
$\sqrt{-75x_1^4-120x_1^3-84x_1^2-28x_1-4}$.
Substituting $x_1=-264/743$, we can check that
\begin{align*}
& \frac{304758098401}{73863713379} f \left(\frac{-264}{743} \right)\\
&+\left(-\frac{115447192511}{73863713379}+\frac{28870417761487643\sqrt{-1}}{27910585919669214} \right) f \left(\frac{-253754+188933 \sqrt{-1}}{863405} \right)\\
& +\left(-\frac{115447192511}{73863713379}-\frac{28870417761487643\sqrt{-1}}{27910585919669214} \right) f \left(\frac{-253754-188933 \sqrt{-1}}{863405} \right)\\
= & -\frac{1}{4 \pi \sqrt{-1}} \int_{S^1} f(z) e^{-\frac{2}{z}} dz 
\end{align*}
is a quadrature formula of degree $4$ for Bessel polynomials with nodes on $\mathbb{Q}(\sqrt{-1})$ corresponding to $P$.
\end{example}

\begin{remark} \label{Rem}
\begin{enumerate}
\item By the same argument as in \cref{MT1}, we can prove that there exists no quadrature formula of degree $2r$ for Bessel polynomials with $r+1$ nodes on $\mathbb{Q}(\sqrt{-11})$ for all $r \in \mathbb{Z}_{\geq 3}$.

\item We cannot prove \cref{MT1} (2) by using %Newton polygons (\cref{NP})
\cref{NP} since there exist quadrature formulas of degree $4$ for Bessel polynomials with $3$ nodes on $\mathbb{Q}(\sqrt{-1})$
by \cref{4-3}.
\end{enumerate}
\end{remark}

\section{Proof of \cref{MT2}}

In this section, we prove \cref{MT2}. 
\begin{proof} [{Proof of \cref{MT2}}]
Note that this is the case $2r+2-k=r+1$. i.e., $k=r+1=3$ in the Riesz--Shohat theorem, \cref{RS}. 
In this case, we prove the following by contradiction: 

There exist no $z_1$, $z_2$, $z_3 \in \mathbb{Q}(\sqrt{-1}) \cap S^1$, $b_1$, $b_2 \in \mathbb{C}$ such that
\[(z-z_1)(z-z_2)(z-z_3)=\phi_3(z)+b_1 \phi_2(z)+ b_2 \phi_1(z).\]
Here, 
\[\phi_n(z):=\frac{y_n(z)}{\frac{(2n)!}{2^n n!}}.\]
Let
\[A:=\begin{pmatrix}
1 & 0 & 0 & 0\\
1 & 1 & 0 & 0\\
\frac{2}{5} & 1 & 1 & 0 \\
\frac{1}{15} & \frac{1}{3} & 1 & 1
\end{pmatrix},\]
and suppose that
\[A\begin{pmatrix}
1\\
b_1 \\
b_2\\
0
\end{pmatrix}
=\begin{pmatrix}
1\\
-z_1-z_2-z_3  \\
z_1z_2+z_2z_3+z_3z_1\\
-z_1z_2z_3
\end{pmatrix}.\]
Comparing the second to the fourth row,
\begin{align*}
b_1 & =-1-z_1-z_2-z_3,\\
b_2 &=\frac{3}{5}+(z_1+z_2+z_3)+(z_1z_2+z_2z_3+z_3z_1),\\
\left(z_1+z_2+\frac{2}{3}+z_1z_2 \right)z_3 &=-\frac{1}{3}-\frac{2}{3} (z_1+z_2)-z_1z_2.
\end{align*}
By the third equality, if
 \[z_1+z_2+\frac{2}{3}+z_1z_2=0,\]
then  
\[z_1z_2=-\frac{2}{3}(z_1+z_2)-\frac{1}{3}.\]
By the above two equalities,
\[z_1+z_2+1=0.\]
Write $z_1=s+t \sqrt{-1}$, $z_2=u+v \sqrt{-1}$ with $s, t, u, v \in \mathbb{Q}$ and $s^2+t^2=u^2+v^2=1$.
Then,  $s+u+1=0$ and $t+v=0$.
Thus, $u=-s-1$ and $v=-t$. Since $s^2+t^2=u^2+v^2=1$,
\[1=(-s-1)^2+t^2=2s+2.\]
Thus, $s=-1/2$ and $t=\pm \sqrt{3}/{2} \not\in \mathbb{Q}$, which is a contradiction.
Therefore,
\[z_1+z_2+\frac{2}{3}+z_1z_2 \neq 0\]
and we have
\begin{align*}
z_3 =&\frac{-\frac{1}{3}-\frac{2}{3} (z_1+z_2)-z_1z_2}{z_1+z_2+\frac{2}{3}+z_1z_2}. 
\end{align*}
We solve $|z_1|^2=|z_2|^2=|z_3|^2=1$ by Mathematica \cite{Wolfram}. 
Here are the outputs: %fn
\fn{$v$ in the first two lines can be expressed as a rational function of $s$, $\sqrt{43+84s+s^2-84s^3-44s^4}$ and $\sqrt{1-s^2}$.} %end

\begin{align*} 
(t,u) &=\left(\sqrt{1-s^2},\frac{-91-202s-112s^2 \pm 2\sqrt{43+84s+s^2-84s^3-44s^4} }{106+224s+120s^2} \right),\\ %\fnm\\
(t,u) &=\left(-\sqrt{1-s^2},\frac{-91-202s-112s^2 \pm 2\sqrt{43+84s+s^2-84s^3-44s^4} }{106+224s+120s^2} \right), \\
(s,t,u,v) &= \left(-1,0, -\frac{1}{2}, \pm \frac{\sqrt{3}}{2} \right),\\
(s,t,u,v) &= \left(1,0, -\frac{9}{10}, \pm \frac{\sqrt{19}}{10} \right),\\
%\end{align*} 
%\begin{align*} 
(s,t,u,v) &= \left(\frac{-28-\sqrt{-11}}{30}, 
-\frac{\sqrt{127-56\sqrt{-11}}}{30}, \frac{-28+\sqrt{-11}}{30}, \frac{\sqrt{127+56\sqrt{-11}}}{30} \right),\\ 
\end{align*} 
\begin{align*} 
(s,t,u,v) &= \left(\frac{-28-\sqrt{-11}}{30}, 
\frac{\sqrt{127-56\sqrt{-11}}}{30}, \frac{-28+\sqrt{-11}}{30},-\frac{\sqrt{127+56\sqrt{-11}}}{30} \right),\\
(s,t,u,v) &= \left(\frac{-28+\sqrt{-11}}{30},
-\frac{\sqrt{127+56\sqrt{-11}}}{30}, \frac{-28-\sqrt{-11}}{30}, \frac{\sqrt{127-56\sqrt{-11}}}{30} \right),\\ 
(s,t,u,v) &= \left(\frac{-28+\sqrt{-11}}{30}, \frac{\sqrt{127+56\sqrt{-11}}}{30}, 
\frac{-28-\sqrt{-11}}{30}, -\frac{\sqrt{127-56\sqrt{-11}}}{30} \right).
\end{align*} 
For the third to the eighth solution, $s \not\in \mathbb{Q}$ or $v \not\in \mathbb{Q}$, which are impossible.

For the first and the second solution, we claim that
\[43+84s+s^2-84s^3-44s^4\]
is a rational square if and only if $s=\pm 1$ or $-1/2$ or $-9/10$, which contradict $t \in \mathbb{Q}$ or $v \in \mathbb{Q}$.

Indeed, let $C$ be an elliptic curve over $\mathbb{Q}$ defined by 
\[y^2=-44x^4-84x^3+x^2+84x+43.\]
By MAGMA \cite{Bosma-Cannon-Playoust}, $\rank C(\mathbb{Q})=0$.
We also know that $(\pm 1,0), \left(-\frac{1}{2}, \pm 3 \right), \left(-\frac{9}{10},\pm \frac{19}{25} \right) \in C(\mathbb{Q})$ 
and the torsion subgroup is isomorphic to $\mathbb{Z}/6 \mathbb{Z}$.
Therefore, we have
\[C(\mathbb{Q})=\left\{(\pm 1,0), \left(-\frac{1}{2}, \pm 3 \right), \left(-\frac{9}{10},\pm \frac{19}{25} \right) \right\}. \]

\end{proof}

\begin{remark} 
\begin{enumerate}
\item
\cref{UQ} implies that there exists a unique quadrature formula of degree $2$ for Bessel polynomials with $2$ nodes on $\mathbb{Q}(\sqrt{-11}) \cap S^1$.
On the other hand, there exists no quadrature formula of degree $3$ for Bessel polynomials with $3$ nodes on $\mathbb{Q}(\sqrt{-11}) \cap S^1$.
The proof is exactly similar to \cref{MT2}.

\item
There exist quadrature formulas of degree $3$ for Bessel polynomials with $3$ nodes on $\mathbb{Q}(\sqrt{-3}) \cap S^1$.
For example,
\begin{align*}
& \frac{2}{3}f(-1)+\frac{3+\sqrt{-3}}{18} f \left(\frac{-1+\sqrt{-3}}{2} \right)+\frac{3-\sqrt{-3}}{18} f \left(\frac{-1-\sqrt{-3}}{2} \right)\\
= &-\frac{1}{4 \pi \sqrt{-1}} \int_{S^1} f(z) e^{-\frac{2}{z}} dz 
\end{align*}
is a quadrature formula of degree $3$ for Bessel polynomials with $3$ nodes on $\mathbb{Q}(\sqrt{-3}) \cap S^1$.
\end{enumerate}
\end{remark}

\section{Proof of \cref{MT3}} 
In this section, we prove \cref{MT3}. 

\begin{proof} [{Proof of \cref{MT3}}]
By \cref{RS}, it is enough to show that there exist $z_1, \ldots, z_{r+1} \in \mathbb{Q}$, $b_1, \ldots, b_r \in \mathbb{C}$ such that
\begin{align} \label{*}
\theta_{r+1}(z):=(z-z_1) \cdots (z-z_{r+1})=\phi_{r+1}(z)+b_1 \phi_r(z) + \cdots +b_r \phi_1(z)+0 \cdot \phi_0(z).
\end{align}
Here, 
\[\phi_n(z):=\frac{y_n(z)}{\frac{(2n)!}{2^n n!}}.\]
Write
\begin{align*} 
(z-z_1) \cdots (z-z_{r+1}) &= 1 \cdot z^{r+1}+Z_1(z_1,\ldots,z_{r+1}) z^r+\cdots +Z_{r+1}(z_1,\ldots,z_{r+1}),\\
\phi_{r+1}(z)+b_1 \phi_r(z) + \cdots +b_r \phi_1(z)&= 1 \cdot z^{r+1}+B_1(b_1,\ldots,b_r) z^r+\cdots +B_{r+1}(b_1,\ldots, b_r).
\end{align*}
Let
\[A:=\begin{pmatrix}
1 & 0 & 0 & 0 & \hdots & 0 & 0 &0\\
a_{(r+1) 1} & 1 & 0 & 0 & \hdots & 0 & 0 &0\\
a_{(r+1) 2} & a_{r1} & 1 & 0 & \hdots & 0 & 0 &0\\
a_{(r+1) 3} & a_{r2} & a_{(r-1)1} & 1 & \hdots & 0 & 0 &0\\
\vdots\ & \vdots & \vdots   & \ddots & \ddots & \vdots & \vdots & \vdots\\
a_{(r+1) (r-1)} & a_{r(r-2)}  & a_{(r-1) (r-3)}  & a_{(r-2) (r-4)} &  
\hdots & 1 & 0 &0\\
a_{(r+1) r} & a_{r(r-1)}  & a_{(r-1) (r-2)}  & a_{(r-2) (r-3)} & \hdots & a_{21} & 1 &0\\
a_{(r+1) (r+1)} & a_{rr}  & a_{(r-1) (r-1)}  & a_{(r-2) (r-2)} & \hdots & a_{22} & a_{11} &1
\end{pmatrix},\]
where we write
\[\phi_k(z):=\sum_{i=0}^k a_{ki} z^{k-i} \quad (a_{k0}=1).\]
Then, by comparing the coefficients of $z^k$ in (2), \cref{*} is equivalent to
\[A\begin{pmatrix}
1\\
b_1 \\
\vdots\\
b_r\\
0
\end{pmatrix}
= \begin{pmatrix}
1\\
B_1(b_1,\ldots,b_r)  \\
\vdots\\
B_{r+1}(b_1,\ldots,b_r) \\
\end{pmatrix} =\begin{pmatrix}
1\\
Z_1(z_1,\ldots,z_{r+1}) \\
\vdots\\
Z_{r+1}(z_1,\ldots,z_{r+1})
\end{pmatrix}.\]
Since $\det(A) \neq 0$, this is equivalent to
\begin{align} \label{**}
\begin{pmatrix}
1\\
b_1 \\
\vdots\\
b_r\\
0
\end{pmatrix}
= A^{-1}\begin{pmatrix}
1\\
Z_1(z_1,\ldots,z_{r+1}) \\
\vdots\\
Z_{r+1}(z_1,\ldots,z_{r+1})
\end{pmatrix}.
\end{align}

The numbers $z_1, \ldots, z_r$,  $b_1, \ldots, b_r$ in \cref{*} can be taken as follows:

First, take $z_1, \ldots, z_r \in \mathbb{Q}$ arbitrarily.
By comparing the last row of \cref{**}, we have a relation among $z_1, \ldots, z_{r+1}$. 
For each $i$, note that $Z_i$ is an elementary symmetric polynomial in $z_1, \ldots, z_{r+1}$ and is linear in $z_{r+1}$.
Since the coefficients of Bessel polynomials are rational, the off-diagonal entries in the lower triangular matrix $A$ are rational numbers.
Thus, the items of the inverse matrix $A^{-1}$ are rational numbers.
Moreover, since the coefficients $Z_1, \cdots, Z_{r+1}$ are linear functions of $z_{r+1}$ with rational coefficients,
 the last row of the equation \cref{**} is a linear equation of $z_{r+1}$ with rational coefficients.
Therefore, $z_{r+1} \in \mathbb{Q}$ can be expressed as a rational function of $z_1, \ldots, z_r$.
Then, by \cref{**}, $b_1, \ldots, b_{r+1}$ can be expressed in terms of $z_1, \ldots, z_{r+1}$.%fn
\fn{Thus, $b_1, \ldots, b_r \in \mathbb{Q}$.} %end
\end{proof}

\begin{example}
The formula
\begin{align*}
 \frac{172}{441} f \left(\frac{1}{2} \right)-\frac{1625}{1611} f \left(\frac{1}{5} \right) + \frac{42592}{26313} f \left(-\frac{27}{44} \right)
= -\frac{1}{4 \pi \sqrt{-1}} \int_{S^1} f(z) e^{-\frac{2}{z}} dz 
\end{align*}
is a quadrature formula of degree $3$ with $3$ nodes on $\mathbb{Q}$.\\
\end{example}

\noindent {\bf Acknowledgements.}
The author thanks anonymous referees for their constructive comments to make this article clearer.
The author thanks his advisor Professor Ken-ichi Bannai for reading the draft and giving helpful comments.
The author also thanks him for warm and constant encouragement.
The author gratefully thanks Professor Yukihiro Uchida for helpful comments and discussions.
The author also thanks Professors Masato Kurihara, Taka-aki Tanaka, Shuji Yamamoto, Yoshinosuke Hirakawa, Kazuki Yamada, Hohto Bekki, Naoto Dainobu and Yoshinori Kanamura for helpful comments and discussions.

\end{document}